\documentclass[12pt]{article}
\usepackage{fullpage}
\usepackage{graphicx}
\usepackage{amsmath,amsthm,amssymb}
\usepackage{color}
\usepackage{url}
\usepackage{hyperref}
\usepackage{euscript}
\usepackage{wasysym}

\newtheorem{theorem}{Theorem}
\newtheorem{lemma}[theorem]{Lemma}

\newtheorem{conj}[theorem]{Conjecture}
\newtheorem{cor}[theorem]{Corollary}

\usepackage{lineno}

\newcommand{\DEFN}[1]{\textit{#1}}

\newcommand\pma{\ensuremath{\mathrm{pm}}}
\newcommand\bpm{\ensuremath{\mathrm{bpm}}}
\newcommand{\MM}{\ensuremath{\EuScript{M}}}

\newcommand{\comment}[1]{}

\title{Deranged matchings: proofs and conjectures}

\begin{document}

 \author{Daniel Johnston\\{\small\texttt{daniel.johnston@trincoll.edu}} \and P. Mark Kayll\\{\small\texttt{mark.kayll@umontana.edu}} \and Cory Palmer\\{\small\texttt{cory.palmer@umontana.edu}}}
 

 
\date{}

\maketitle

\begin{abstract}
\noindent
  We introduce, and partially resolve, a conjecture that brings a
  three-centuries-old derangements phenomenon and its much younger
  two-decades-old analogue under the same umbrella. Through a
  graph-theoretic lens, a derangement is a perfect matching in the
  complete bipartite graph $K_{n,n}$ with a disjoint perfect matching
  $M$ removed. Likewise, a deranged matching is a perfect matching in
  the complete graph $K_{2n}$ minus a perfect matching $M'$. With
  $\pma(\cdot)$ counting perfect matchings, the elder phenomenon
  takes the form $\pma(K_{n,n}-M)/\pma(K_{n,n})\to 1/e$ as
  $n\to\infty$ while its youthful analogue is
  $\pma(K_{2n}-M')/\pma(K_{2n})\to 1/\sqrt{e}$. These starting graphs
  are both $2n$-vertex `balanced complete $r$-partite' graphs 
  $K_{r \times {2n}/{r}}$, respectively with $r=2$ and $r=2n$. We
  conjecture that
  $\pma(K_{r\times{2n}/r}-M)/\pma(K_{r\times{2n}/r})\sim e^{-r/(2r-2)}$ as
  $n\to\infty$ and establish several substantive special cases
  thereof. For just two examples, $r=3$ yields the limit $e^{-3/4}$ while $r=n$
  results again in $e^{-1/2}$.  Our tools blend combinatorics and
  analysis in a medley incorporating Inclusion-Exclusion
  and Tannery's Theorem.
\end{abstract}

\begin{quote}
Such permutations are often called \textit{derangements}, and counting
them is a traditional preoccupation of Combinatorics texts.
\hfill\textsc{C.D.~Godsil},~1993~\cite{Godsil93}
\end{quote}

\renewcommand{\thefootnote}{}
\footnotetext{\emph{Date:} 31 October 2022 \raisebox{-2pt}{{\Large $\blacksmiley$}}}
\footnotetext{Partially supported by grants from the Simons Foundation (\#279367 to Mark Kayll and \#712036 to Cory Palmer).}
\footnotetext{Copyright \copyright\ 2022 by the authors}
\renewcommand{\thefootnote}{\arabic{footnote}}

\section{Introduction}

The functions in the opening quote are permutations of the set
$\{1,\ldots,n\}$ having no fixed points. They are perhaps the first
nontrivial objects ever counted, having been introduced by R\'{e}mond
de Montmort in his famous 1708 treatise~\cite{Montmort1708-13} and
determined in its 1713 second edition.
Incidentally, Sir Isaac Newton owned a copy of the latter work, a
print from which appears as the frontispiece of Bollob\'{a}s'
monograph~\cite{Bollobas86}.  So by now it might be surprising to
learn that anybody has anything new to share about derangements.

`Everybody' knows the beautiful appearance of Euler's number $e$ in their
enumeration. As a reminder, the number $d_n$ of derangements of
$\{1,\ldots,n\}$ is the closest integer to $n!/e$; equivalently,
\begin{equation}
  \label{derange-formula-basic}
d_n = n! \sum_{k=0}^n \frac{(-1)^k}{k!}.
\end{equation}
These relations imply that of \emph{all}
the permutations of $\{1,\ldots,n\}$, about a fraction $1/e$ of them
are derangements. 

It's fair to say that most mathematics students see the
entrance of $e$ into the discussion of $d_n$ as unexpected, but it is
actually quite natural. If we think of assembling a permutation $\pi$ of
$\{1,\ldots,n\}$ by letting each $i$ in this domain decide its image,
then for the resulting $\pi$ to be a derangement, each $i$ has 
$n-1$ choices for $\pi(i)$ because it can't map to itself. If the $n$
decisions had no interdependence (hold it---we know that's too big an
``If''), then the overall proportion of assembled permutations that are
deranged would be
\[
\left(\frac{n-1}{n}\right)^{n}=\left(1-\frac{1}{n}\right)^{n}\to\frac{1}{e} ~~\text{(as $n\to\infty$)}.
\]
Of course, the requirement that $\pi$ be surjective---in addition merely
to moving each domain element $i$---enters dependence into those $n$
decisions. But, especially when $n$ is large, the error in the heuristic
thinking is not large, so rigour can be salvaged. Indeed, other recent
\textsc{Monthly} authors~\cite{Bary-SorokerGorodetsky2018} have
accomplished this exact feat in an enlightening analysis of proximate
probability measures. Though we shall not pursue this angle, the
heuristic did play an essential role in formulating 
Conjecture~\ref{kayll-conj} below.

We shall approach derangements through a graph-theoretic lens. First
notice that each permutation $\pi$ of $\{1,\ldots,n\}$ can be viewed as a
perfect matching in the complete bipartite graph $K_{n,n}$ (with each
$i$ in the domain matched to its image $\pi(i)$ across the
bipartition). If $\pi$ happens to be a derangement, then, as a perfect
matching, it includes none of the edges of the particular perfect matching
$M=\{11',22',\ldots,nn'\}$ in $K_{n,n}$ and thus is a perfect matching
in $K_{n,n}-M$. Conversely, every perfect matching in the latter graph
arises from a derangement.

Following Godsil~\cite{Godsil93}, denote the number of perfect
matchings in a graph $G$ by $\pma(G)$. 
The preceding remarks give $\pma(K_{n,n})=n!$ and
$\pma(K_{n,n}-M)=d_n$, and since all perfect matchings of $K_{n,n}$
are isomorphic, the $M$ in this last identity can be any fixed perfect
matching of this graph. These observations, together with
(\ref{derange-formula-basic}), show that
\begin{equation}
  \label{hat-check-graphs}
\frac{\pma(K_{n,n}-M)}{\pma(K_{n,n})} =\frac{d_n}{n!}\to e^{-1}~~\text{(as $n\to\infty$)}.
\end{equation}

Determining $d_n$ is the `Hat-Check Problem' (see, e.g.,
\cite{West2021}), which, as noted in the first paragraph, is
more than 300 years old. Its advanced age heightens the surprise when
one learns---perhaps right now---that it took almost this long for
the following analogue to attract study.

\begin{center}
\textbf{Kindergartner Problem}
\begin{quote}
  An even number $2n$ of kindergartners buddy up for a
  field trip; at the end of the day, their inexperienced chaperone
  asks them to buddy up again arbitrarily. What's the chance that no
  original buddy pairs persist?
\end{quote}
\end{center}

\noindent
If the children are the vertices of the complete graph $K_{2n}$, then
the original pairing is a perfect matching $M$ in $K_{2n}$, and the
problem asks about the ratio of $\pma(K_{2n}-M)$ to
$\pma(K_{2n})$. With $(\cdot)!!$ denoting the double factorial
function, the latter count is $\pma(K_{2n})=(2n-1)!!$.
The matchings counted by $\pma(K_{2n}-M)$ are called
\textit{deranged matchings}; for comparison with
(\ref{hat-check-graphs}), we also use $D_n$ to denote their
count. Brawner~\cite{Brawner2000} conjectured the following analogue
of (\ref{hat-check-graphs}):
\begin{equation}
  \label{kinder}
\frac{\pma(K_{2n}-M)}{\pma(K_{2n})} =\frac{D_n}{(2n-1)!!}\to e^{-1/2}~~\text{(as $n\to\infty$)}.
\end{equation}

Published in 2000, this conjecture arose from an inquiry by the United
States Tennis Association concerning the tournament draw for the 1996
U.S. Open. The kindergartner-formulation is due to the second author
of the present article and first appeared in \cite{Morey-thesis-2013}; see
also \cite{Morey2013}. The conjecture was proved by
Margolius~\cite{Margolius2001} and again by the second author in this
\textsc{Monthly}~\cite{Kayll-AMM}.

In their original formulations, these problems weren't
framed to display their similarity so transparently as in
(\ref{hat-check-graphs}) and (\ref{kinder}). Now that they are, hope
springs to put them under one umbrella.

Besides both containing $2n$ vertices, the graphs $K_{n,n}$ and
$K_{2n}$ share a deeper connection, namely being examples of 
`balanced complete $r$-partite' graphs. We spell this out. 
For a positive integer $r$, an \DEFN{$r$-partite} graph $G$ is one whose
vertex set $V$ can be partitioned into $r$ classes,
$V=V_1\cup\cdots\cup V_r$, so that no edge of $G$ has both ends in the
same partition class. If, moreover, $G$ contains every conceivable edge $xy$
with $x\in V_i$ and $y\in V_j$ (whenever $i\neq j$), then $G$ is
\DEFN{complete $r$-partite}. Finally, such a graph is \DEFN{balanced} when each
class $V_i$ contains the same number of vertices. 
Following \cite{BrouwerCohenNeumaier89}, we use 
$K_{r \times {2n}/{r}}$ to denote the $2n$-vertex balanced complete
$r$-partite graph. Thus, $K_{n,n}=K_{2\times n}$ and $K_{2n}=K_{2n\times 1}$.

The latter observation and the enticing similarity between the phenomena
(\ref{hat-check-graphs}) and (\ref{kinder}) led the second author to
apply the heuristic sketched above to general balanced complete
$r$-partite graphs. This resulted in the following `linking
conjecture', in which we use standard asymptotic
tilde-notation.\footnote{Real sequences $(a_n)$ and $(b_n)$ satisfy 
$a_n\sim b_n$ (as $n\to\infty$) exactly when $\lim_{n\to\infty}(a_n/b_n)=1$.}
\begin{conj}[\cite{Kayll2014-Fulbright-prop}]
\label{kayll-conj}
If $r=r(n)\geq 2$ is integer-valued and divides $2n$, and
$M$ is a perfect matching in $K_{r \times {2n}/{r}}$, then
\[
 \frac{\pma(K_{r \times {2n}/{r}}-M)}{\pma(K_{r \times {2n}/{r}})} \sim e^{-r/(2r-2)}~~\text{(as $n\to\infty$)}.
\]
\end{conj}

The solutions to the Hat-Check and Kindergartner
Problems confirm Conjecture~\ref{kayll-conj} for $r=2$ and $r=2n$, respectively. In
this article, we give further supporting evidence.
Section~\ref{const-classes-section} verifies Conjecture~\ref{kayll-conj} 
for $r=3$ and shows that it holds for constant $r\geq 4$ when restricted to
certain `nice' perfect matchings $M$.
Section~\ref{linear-classes-section} settles the conjecture
when $r(n)$ is linear in $n$, i.e., when the classes of
$K_{r \times {2n}/{r}}$ are of constant size. Up to, but not including, 
Section~\ref{refinement-section}, our presentation is self-contained. In
that short section, we show how a deep result of
McLeod~\cite{McLeod2009} allows us to extend the resolution of  
Conjecture~\ref{kayll-conj} to cases of $r(n)$ dropping to fractional
powers $\Omega(n^{\delta})$ (for fixed positive $\delta<1$). Thus, we leave
open essentially just the cases of polylogarithmic $r(n)$ (beyond the
general constant $r\geq 4$ noted above).

\subsection*{Background and methods}

As noted in the introductory paragraph, the study of derangements is
generally believed to have been initiated by
Montmort~\cite{Montmort1708-13}, who in 1708
described a probabilistic game of coincidences.
He analyzed his ``Jeu du Treize'' by an iterative procedure in 1713,
when he also included a
solution by his friend Nicolas Bernoulli [1687--1759] (we
mention his lifespan to identify the member
of this extended mathematical family).  Bernoulli's
solution---addressed more recently in \cite{Dupont78}---is
a classical application of the Principle of Inclusion-Exclusion,
which appears as Theorem~\ref{PIE-stmt} below.

Other approaches leading to (\ref{derange-formula-basic}) abound.  For
example, one can deduce combinatorially the recurrence relation
$d_n = (n-1)(d_{n-1}+d_{n-2})$ (as did Euler), then solve it iteratively, or by
substitution, or by incorporating generating functions. Likewise for
the recurrence $d_n = nd_{n-1}+(-1)^n$.  Another proof
applies the Binomial Inversion Theorem  
to the identity $n! = \sum_{k=0}^n \binom{n}{k} d_k$
(see, e.g., \cite[Chapter~5]{GKP}). Two others derive from the
relation $d_n=\int_{0}^{\infty}(t-1)^ne^{-t}dt$, which follows from
a result of Joni and Rota~\cite{JoniRota80} and Godsil~\cite{Godsil81}
(because $(t-1)^n$ is the `rook polynomial' of a single perfect
matching of $n$ edges in an otherwise empty graph---see
\cite{Godsil93} or, e.g.,
\cite{EmersonKayll} or \cite{Kayll2011}). An analogous 
integral formula for counting perfect matchings in general graphs
provides the basis for Lemma~\ref{godsil-lemma} below.

Though we're not attempting to be encyclopedic, the preceding
paragraphs give a reasonably complete overview of the basic techniques
leading to (\ref{derange-formula-basic}). For more on the math, we
refer the reader to the references already cited and to
\cite{Stanley-v1-2012} or \cite{West2021}.
Derangements fall in the realm of `permutations with restricted
position', a branch of combinatorics originally exposited by
Riordan~\cite{Riordan58}; this reference also provides some related
history. We found the thesis~\cite{Bhatnagar95-thesis} and
article~\cite{HansonSeyffarthWeston83} similarly helpful in this regard.

The latter article presents a generalization of the Hat-Check Problem
different from ours. Others have done this, too. David~\cite{David88}
views the problem through the lens of counting certain elements of the
symmetric group on $n$ symbols. In his generalization and treatment,
the Poisson probability distribution---and hence powers of $e$ not
unlike ours---appear(s).  Penrice~\cite{Penrice91} considers a ground
set $S$ partitioned into equal-sized parts and determines the limiting
probability that a random permutation of $S$ maps every element of $S$
outside its partition class. Powers of $e$ arise in this generalization as well.

Clark~\cite{Clark2013} introduces a generalization with quite a
different character.  `Living' on a simple graph $G=(V,E)$, his
\DEFN{graph derangements} are bijections $f$ of $V$ to itself for
which each $v\in V$ is adjacent to its image $f(v)$.
\DEFN{Graph permutations} are defined likewise, except vertices are
also allowed to remain fixed.  Clark recovers ordinary derangements
(resp.\ permutations) when $G=K_n$. Aside from establishing some basic
properties of this version of derangement, his definition inspired
at least a couple of recent studies: \cite{AusthofBennettChristo2022} and
\cite{BDHHL2021}. These explore the ratios of
derangement- to permutation-counts in various graphs and digraphs.

Like the published literature, the Internet is replete with
information about derangements. For example, the websites
\cite{EncyclopediaMath-page} and \cite{MacTutor-page} provide detailed
historical accounts of Montmort. Not surprisingly, there are also
separate \textsf{Wikipedia}\textregistered\ pages on Montmort and derangements.  And the
numbers $d_n$, $D_n$, respectively, of derangements and deranged
matchings appear as sequences \textsf{A000166}~\cite{OEIS-A000166},
\textsf{A053871}~\cite{OEIS-A053871} in The On-Line Encyclopedia of
Integer Sequences\textregistered.

One can observe in many of these earlier sources a shared theme but
differing substance from the present paper.  These feel tantalizingly
close to our work here, but mathematically they stand separately.

The Principle of Inclusion-Exclusion (PIE) furnishes our primary
combinatorial tool while Tannery's Theorem proves fruitful on the
analysis front. We discuss these results in Section~\ref{sect-prelims}, but
see \cite{West2021} and \cite{aaaa}, respectively, for thorough
treatments. For any omitted graph theory terminology,
almost any source will suffice (e.g., \cite{BondyMurty2008}).

\section{Preliminaries}
\label{sect-prelims}

This section records the main tools just mentioned, starting 
with PIE. Of the various ways to formulate this principle, the most convenient here is to count the elements in the complement of a finite union of finite sets in terms of those sets' various intersections. So the initial data includes a finite index set $I$ and, for each $t \in I$, a subset $A_t$ of a finite `universal' set $U$. With this backdrop, we then have the following 

\begin{theorem}[Principle of Inclusion-Exclusion]
\label{PIE-stmt}
The number of elements of $U$ not contained in any $A_t$ is given by

\hfill\(\displaystyle{
\left|\left(\bigcup_{t\in I}A_t \right)^C\right| =
\sum_{S \subseteq I} (-1)^{|S|} \left| \bigcap_{t \in S} A_t \right|.
}\)\hfill$\Box$
\end{theorem}

Our other main tool is a useful, though not widely known, fact from
analysis. It's a special case of a famous result---Lebesgue's Dominated Convergence Theorem (see, e.g., \cite{aaaa})---but, especially for those mathematicians working in the 
`discrete world', it deserves its own limelight.
For the sake of completeness, and maybe a bit of advertising, we include a short proof. The result concerns the legality of interchanging the order of the limit and summation operations
when the terms themselves are not assumed constants but are converging. Here the initial data includes a sequence
$(f_i(\cdot))_{i\geq 1}$ of real-valued functions on 
$\mathbb{N}$ together with two real-valued sequences $(f_i)$, $(M_i)$ (the latter serving, respectively, as limits and absolute upper bounds for the evaluation sequences 
$(f_{i}(n))_{n\geq 1}$).

\begin{lemma}[Tannery's Theorem]\label{analysis-lemma-1}
    If, for each $i\geq 1$, we have
    $\lim_{n\to\infty}f_{i}(n)=f_{i}$ and 
    $|f_i(n)| \leq M_i$ for all $n \in \mathbb{N}$, and
	furthermore if $\sum_{i=1}^\infty M_i < \infty$,
	then for any integer-valued sequence $(s_n)$ with 
	$s_n\to\infty$, we have
	\begin{equation}
	\label{eq-tannery-conclusion}
	\lim_{n \rightarrow \infty} \sum_{i=1}^{s_n} f_{i}(n) = \sum_{i=1}^\infty  f_{i}.
	\end{equation}
\end{lemma}

\noindent
We choose the conclusion's version in 
(\ref{eq-tannery-conclusion}) because,
among several alternatives, it's the most natural for our
applications. 

\begin{proof} 
Observe that $|f_i(n)| \leq M_i$ implies that $|f_i| \leq M_i$.
Therefore, by the comparison test, $\sum |f_i|$ converges, and so must $\sum f_i$. 

	Fix $\epsilon > 0$. The convergence of $\sum M_i$ implies that there is an integer $N= N(\epsilon)$ such that 
		\[
	 \sum_{i = N+1}^\infty M_i  < \frac{\epsilon}{4}.
	\]
	Now pick $N'=N'(\epsilon)$ such that 
	whenever $n > N'$, we have $s_n > N$ and
	\[
	|f_i(n) - f_i| < \frac{\epsilon}{2N}
	\]
for every $i\in\{1,2,\dots, N\}$.	

	Therefore, if $n > N'$, then
	\begin{align*}
	\left|\sum_{i=1}^{s_n} f_{i}(n) - \sum_{i = 1}^\infty f_{i} \right| 
	& = \left| \sum_{i=1}^{N} (f_{i}(n)- f_i) +  \sum_{i=N+1}^{s_n}  f_{i}(n) - \sum_{i=N+1}^{\infty}   f_{i}\right|    \\
	& \leq \sum_{i=1}^{N} \left|f_{i}(n)- f_i\right| +   \sum_{i=N+1}^{s_n}  |f_{i}(n)| +\sum_{i=N+1}^{\infty}   |f_{i}|    \\
	& <  N \frac{\epsilon}{2N}  + 2 \sum_{i=N+1}^\infty M_i \\
	& < \frac{\epsilon}{2} + \frac{\epsilon}{2} \\
	& =  \epsilon. \qedhere\\
	\end{align*}    
\end{proof}

In the proofs of 
Theorems~\ref{3-part-thm} and \ref{k-part-thm}
(Section~\ref{const-classes-section}), we shall need a 
slightly more general version of Tannery's Theorem: one that disregards the order of the summation index. We close this section 
with a statement and proof of the generalization. To this end,
let us recall that for a monotone increasing sequence $(X_n)$ of sets,  we have $\lim_{n\to\infty} X_n = \bigcup_{n=1}^{\infty} X_n$. Furthermore, for a countable set $X$, the sum $\sum_{x \in X} f(x)$ is defined only when there is a bijection $g\colon\mathbb{N} \to X$ such that the series $\sum_{n=1}^\infty f(g(n))$ converges absolutely.

\begin{cor}
[Tannery's Theorem -- unordered index version]\label{analysis-lemma-2}
    Let $I_1 \subset I_2 \subset \cdots$ be a monotone increasing sequence of finite subsets of a countable set $I$ such that $\lim_{n\to\infty} I_n = I$.
	If, for each $t \in I$, we have
	$\lim_{n\to\infty}f_{t}(n)=f_{t}$ and 
	$|f_t(n)|\leq M_t$ for all $n\in\mathbb{N}$, and
	furthermore if $\sum_{t\in I}M_t<\infty$, then
	\[
	\lim_{n \to\infty}\sum_{t\in I_n}f_{t}(n)=\sum_{t\in I}f_{t}.
	\]
\end{cor}

\begin{proof}
By the hypotheses on $(I_n)$, there is an enumeration of $I$ respecting the monotonicity of $(I_n)$, i.e., a bijection
 $g \colon \mathbb{N} \to I$ such that
 $\{g(1),g(2),\dots, g(|I_n|)\} = I_n$ for all $n \in \mathbb{N}$. 
Consequently, each $f_{g(i)}(n)$ converges to $f_{g(i)}$ as $n \rightarrow \infty$, each $|f_{g(i)}(n)| \leq M_{g(i)}$ for all $n \in \mathbb{N}$, and $\sum M_{g(i)} < \infty$.  Therefore, we may apply Lemma~\ref{analysis-lemma-1} with $s_n=
|I_n|$ to obtain
\[
    \lim_{n \rightarrow \infty} \sum_{t \in I_n} f_t(n) = \lim_{n \rightarrow \infty} \sum_{i=1}^{|I_n|} f_{g(i)}(n) = \sum_{i=1}^\infty f_{g(i)} = \sum_{t \in I} f_t. \qedhere \\
\]
\end{proof}

\section{Instances of constant {\boldmath $r(n)$}} 
\label{const-classes-section}

The proof of Conjecture~\ref{kayll-conj} when $r=r(n)=3$ sheds 
lots of light on the cases when $r$ takes any constant value 
(dividing $2n$), but shadows remain. In this section, we first present that proof (Theorem~\ref{3-part-thm}) and then indicate 
how it lifts to certain other cases of constant $r$ 
(Theorem~\ref{k-part-thm}).

The hypothesis of $r$ dividing $2n$ when $r=3$ implies that 
$n$ is restricted to values such that $2n=6m$ for some positive integer $m$. In this case, the balanced complete $r$-partite
graph $K_{r\times 2n/r}=K_{3\times 6m/3}$ is the complete 
tripartite graph $K_{2m,2m,2m}$. As this notation is not 
yet too unwieldy, we adopt it for establishing this case of 
Conjecture~\ref{kayll-conj}. (But we draw the line before 
turning to 
$K_{\underbrace{\scriptstyle{(r-1)m,(r-1)m,\ldots,(r-1)m}}_r}$.)

\begin{theorem}\label{3-part-thm}
For any given perfect matching $M$ in the complete tripartite graph $K_{2m,2m,2m}$, we have
\[	
\lim_{m \rightarrow \infty} \frac{\pma(K_{2m,2m,2m}-M)}{\pma(K_{2m,2m,2m})} = e^{-3/4}.
\]
\end{theorem}

\begin{proof}
    Let us begin by counting the number of perfect matchings in
    $K_{2m,2m,2m}$, whose partition classes we'll denote by 
    $X,Y,Z$.
    Observe that $M$ has $3m$ total edges and each partition class is incident to exactly $2m$ of these edges (as illustrated in Figure~\ref{first-figure}).
    This implies that in $M$, there are exactly $m$ edges between each pair of classes. We  imagine building $M$ step-by-step.
Its edges effectively partition the $2m$ vertices of $Y$ into two
sets, one of size $m$ (matched with vertices of $X$) and one of size
$m$ (matched with vertices of $Z$).  This can be done in
$\binom{2m}{m}$ ways. The classes $X$ and $Z$ are partitioned in a
similar fashion. With these three partitions set, there are $m!$ ways
to match the vertices of $Y$ and $X$, $m!$ ways to match those of $Y$
and $Z$, and $m!$ ways to match those of $X$ and $Z$ to finish
building $M$. 
Thus,
    \begin{equation}\label{PM-in-tripartite}
        \pma(K_{2m,2m,2m}) = \binom{2m}{m}^3(m!)^3.
    \end{equation}
    
\begin{figure}[ht]
    \centering
    	\includegraphics[scale=.3]{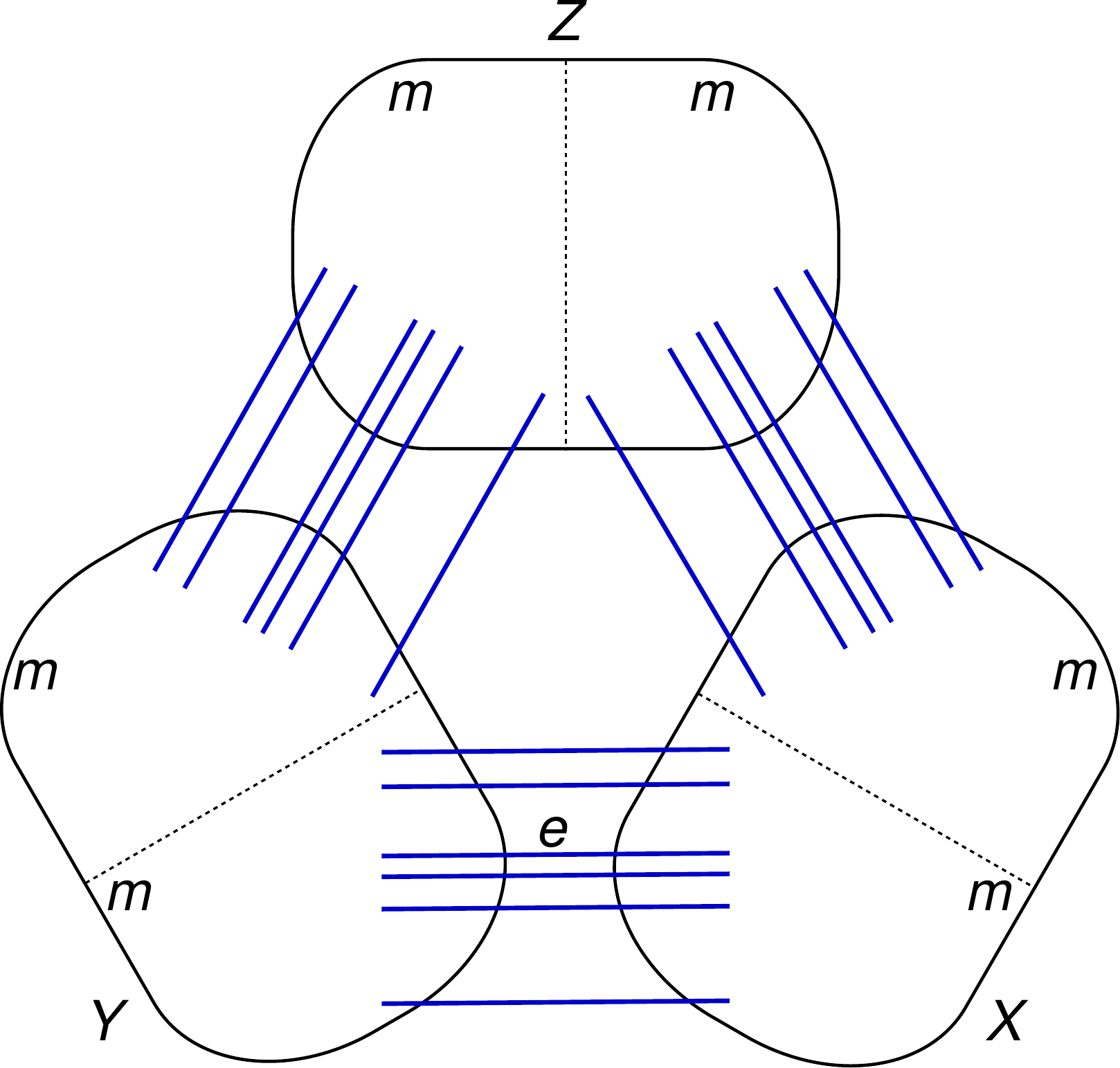}
    \caption{A perfect matching $M$ in $K_{2m,2m,2m}$}
    \label{first-figure}
\end{figure}
    
	Now fixing a particular $M$, we'll use inclusion-exclusion (Theorem~\ref{PIE-stmt}) to
	count the number of perfect matchings in $K_{2m,2m,2m}-M$. 
	For an edge $e \in M$, let $A_e$ be the family of perfect matchings in $K_{2m,2m,2m}$ that use edge $e$. 
	Thus, for a subset $L \subseteq M$, the intersection
	$\bigcap_{e \in L} A_e$ is the family of perfect matchings in $K_{2m,2m,2m}$ that include all of the edges of $L$.
    By Theorem~\ref{PIE-stmt}, the number of perfect matchings in $K_{2m,2m,2m}$ that do not use any edge of $M$ is
	\[
	\pma(K_{2m,2m,2m}-M) = \sum_{L \subseteq M} (-1)^{|L|} \left| \bigcap_{e \in L} A_e \right|.
	\]
	
    It's helpful to expand this sum by specifying the edges  of $M$ between pairs of  partition classes.
    Write $M=M_{XY} \cup M_{YZ} \cup M_{XZ}$, where $M_{AB}$ denotes those $M$-edges between classes $A$ and $B$. 
    Thus we have
    \begin{equation}\label{PIE}
     \pma(K_{2m,2m,2m}-M) =\sum_{I \subseteq M_{XY}} \sum_{J \subseteq M_{YZ}} \sum_{K \subseteq M_{XZ}} (-1)^{| I \cup J \cup K|} \left| \bigcap_{e \in I} A_e \cap \bigcap_{f \in J} A_f \cap \bigcap_{g \in K} A_g \right|.
    \end{equation}
    A term in this sum counts the number of perfect matchings in $K_{2m,2m,2m}$ that include all edges in 
    $L=I \cup J \cup K$; see Figure~\ref{second-figure}, where the ellipses capture $L$.

    \begin{figure}[ht]
    \centering
    	\includegraphics[scale=.3]{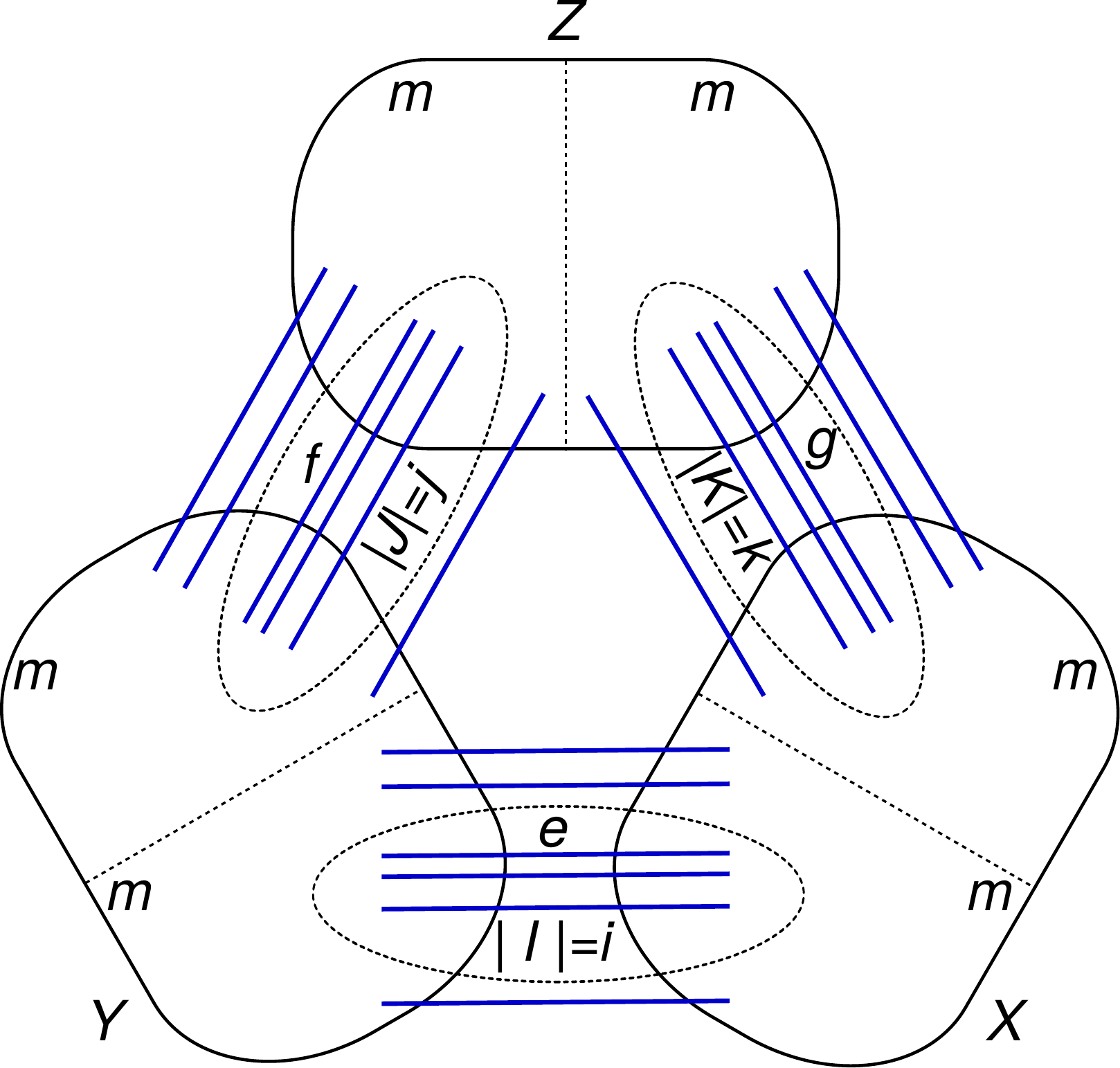}
    \caption{Highlighted edges $I\cup J\cup K$ of the deleted  perfect matching}\label{second-figure}
\end{figure}
    
    To evaluate these terms, we count the number of ways to build a perfect matching $M'$ that includes the edges of $I \cup J \cup K$. Figure~\ref{third-figure} depicts the various parameters in the following discussion.
    First consider the class $Y$, which has $|I| + |J|$ vertices already matched by edges of $I \cup J$.  
    Since $M'$ must have exactly $m$ edges between each pair of classes, its edges partition the remaining unmatched 
    $2m - |I| - |J|$ vertices of $Y$ into two sets, one of size 
    $m - |I|$ (matched with vertices of $X$) and one of size 
    $m - |J|$ (matched with vertices of $Z$).  This can be done in $\binom{2m-|I|-|J|}{m-|I|}$ ways. The classes $Z$ and $X$ are partitioned in a similar fashion. Once these partitions are set, there are $(m-|I|)!$ ways to match the  as-yet unmatched vertices of $Y$ and $X$, $(m-|J|)!$ ways to match those vertices of $Z$ and $Y$, and $(m-|K|)!$ ways to match those vertices of $X$ and $Z$ to finish building $M'$. Therefore, a term $\left| \bigcap_{e \in I} A_e \cap \bigcap_{f \in J} A_f \cap \bigcap_{g \in K} A_g \right|$ is equal to
    \[
    \binom{2m-|I|-|J|}{m-|I|} \binom{2m-|J|-|K|}{m-|J|}\binom{2m-|K|-|I|}{m-|K|} (m-|I|)!(m-|J|)!(m-|K|)!.
    \]
    
    \begin{figure}[ht]
    \centering
    	\includegraphics[scale=.3]{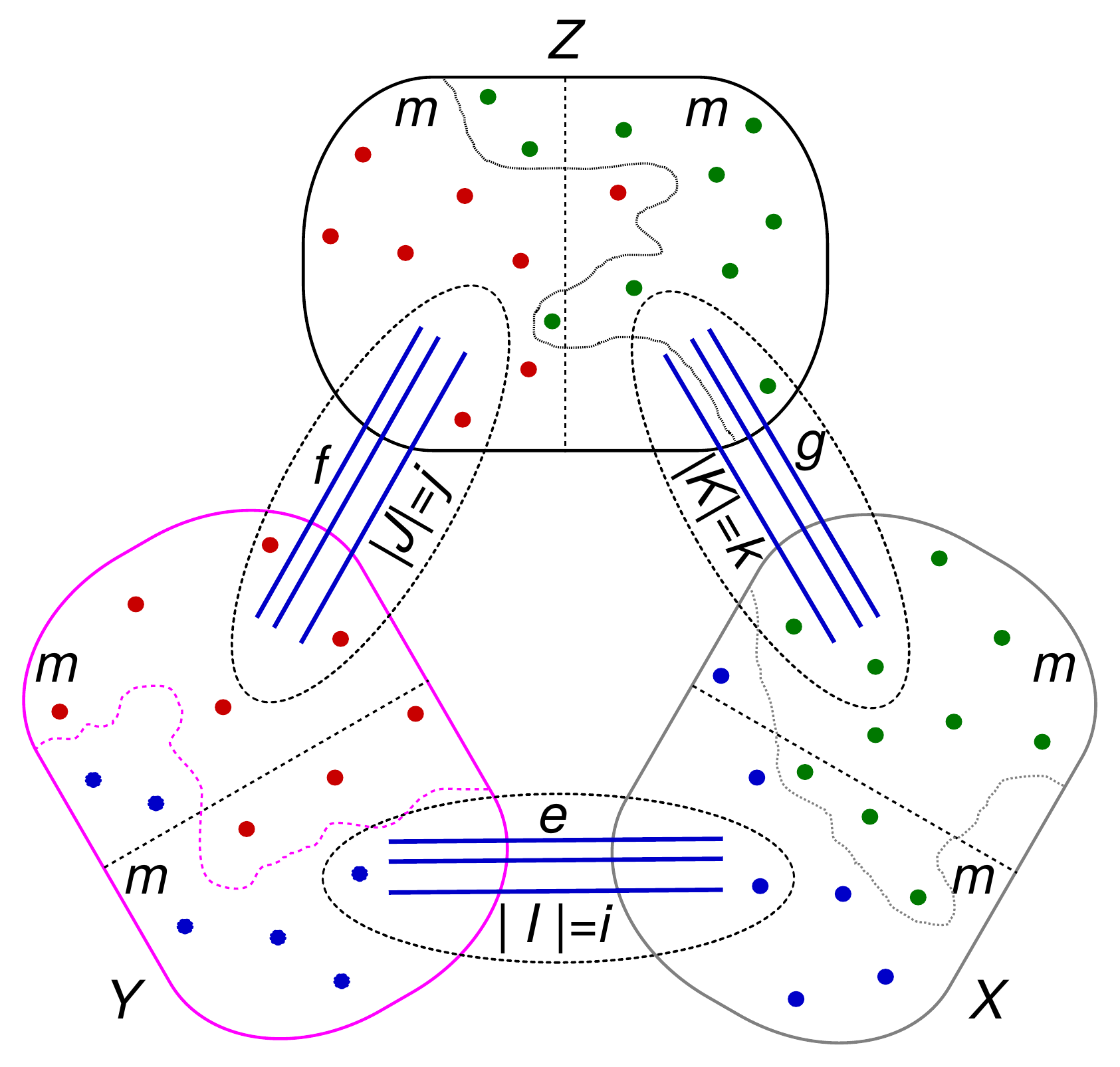}
    \caption{Building a perfect matching in $K_{2m,2m,2m}-M$}\label{third-figure}
\end{figure}

   Observe that the expression above depends only on the cardinalities of $I$, $J$, and $K$. Therefore, we can reindex the sum on the right side of (\ref{PIE}) to express $\pma(K_{2m,2m,2m}-M)$ as
\begin{equation}\label{run-on-equation}
\begin{split}
          \sum_{0 \leq i,j,k\leq m} (-1)^{i+j+k} \binom{m}{i} \binom{m}{j}\binom{m}{k}\binom{2m-i-j}{m-i}&\binom{2m-j-k}{m-j}\binom{2m-k-i}{m-k} \times \\
      & \qquad \qquad  \times (m-i)!(m-j)!(m-k)!.
\end{split}
\end{equation}

\noindent
Dividing this sum by
$\pma(K_{2m,2m,2m}) = \binom{2m}{m}^3(m!)^3$ and simplifying gives
\begin{align*}
\frac{\pma(K_{2m,2m,2m}-M)}{\pma(K_{2m,2m,2m})} &=  \sum_{0 \leq i,j,k \leq m} \frac{(-1)^{i+j+k}}{i!j!k!}\frac{\binom{2m-i-j}{m-i}}{\binom{2m}{m}} \frac{\binom{2m-j-k}{m-j}}{\binom{2m}{m}}\frac{\binom{2m-k-i}{m-k}}{\binom{2m}{m}}.
\end{align*}

To determine the limiting value of this ratio, we aim to apply
Corollary~\ref{analysis-lemma-2}. Any pair $a,b \in \{i,j,k\}$
satisfies
\[
\lim_{m \rightarrow \infty} \frac{\binom{2m-a-b}{m-a}}{\binom{2m}{m}} = \lim_{m \rightarrow \infty} \frac{\prod_{\ell=0}^{a-1} (1-\ell/m) \prod_{\ell=0}^{b-1} (1-\ell/m)}{\prod_{\ell=0}^{a+b-1} (2-\ell/m)}=  \frac{1}{2^{a+b}}.
\]
This implies that for all $i,j,k$, we have
\begin{equation}\label{3-part-term}
 \lim_{m \to \infty} \frac{(-1)^{i+j+k}}{i!j!k!}\frac{\binom{2m-i-j}{m-i}}{\binom{2m}{m}} \frac{\binom{2m-j-k}{m-j}}{\binom{2m}{m}}\frac{\binom{2m-k-i}{m-k}}{\binom{2m}{m}} = \frac{(-1)^{i+j+k}}{i!j!k!} \frac{1}{2^{2i+2j+2k}} = \frac{(-1)^{i+j+k}}{i!j!k!} \frac{1}{4^{i+j+k}}.
\end{equation}

\noindent
Moreover, the absolute value of the expression inside the limit
in (\ref{3-part-term}) is at most $\frac{1}{i!j!k!}$, and 
\[
\sum_{0 \leq i,j,k \leq \infty} \frac{1}{i!j!k!} =\lim_{m \rightarrow \infty} \sum_{0 \leq i,j,k \leq m} \frac{1}{i!j!k!} = \lim_{m \rightarrow \infty} \left( \sum_{i=0}^m \frac{1}{i!} \sum_{j=0}^m \frac{1}{j!} \sum_{k=0}^m \frac{1}{k!} \right)= e^3 < \infty.
\]

\medskip
The  three assertions starting  at  (\ref{3-part-term}) establish the hypotheses of Corollary~\ref{analysis-lemma-2},  and therefore
\begin{align*}
\lim_{m \rightarrow \infty} \frac{\pma(K_{2m,2m,2m}-M)}{\pma(K_{2m,2m,2m})} & = \nonumber
\lim_{m \rightarrow \infty} \sum_{0 \leq i,j,k \leq m} \frac{(-1)^{i+j+k}}{i!j!k!}\frac{\binom{2m-i-j}{m-i}}{\binom{2m}{m}} \frac{\binom{2m-j-k}{m-j}}{\binom{2m}{m}}\frac{\binom{2m-k-i}{m-k}}{\binom{2m}{m}} \\ \nonumber
& = \lim_{m \rightarrow \infty} \sum_{0 \leq i,j,k \leq m} \frac{(-1)^{i+j+k}}{i!j!k!} \frac{1}{4^{i+j+k}} \\
& = \lim_{m \to \infty}  \left( \sum_{i=0}^m \frac{(-1)^{i}}{i!4^i}\sum_{j=0}^m  \frac{(-1)^{j}}{j!4^j}\sum_{k=0}^m \frac{(-1)^{k}}{k!4^k} \right) \\ \nonumber
&= e^{-1/4}e^{-1/4}e^{-1/4} \\ \nonumber
& = e^{-3/4}. \qedhere \nonumber
\end{align*}
\end{proof}

As noted at the start of the preceding proof,
a perfect matching in a balanced complete tripartite graph always has the same number of edges between each pair of classes. Unfortunately, this property is not shared by a balanced complete $r$-partite graph for $r \geq 4$. 
Indeed, consider the complete $4$-partite graph with classes $V_1,V_2,V_3,V_4$, each of size $m$. For example, there is a perfect matching with $m$ edges between $V_1$ and $V_2$ and $m$ between $V_3$ and $V_4$ and consequently no edges between 
$V_1$ and $V_3$ (nor between $V_2$ and $V_4$, etc.).

We've arrived at the edge of our Conjecture~\ref{kayll-conj} shadows. In order to generalize 
Theorem~\ref{3-part-thm}---and hence confirm
Conjecture~\ref{kayll-conj} for other constants 
$r\geq 4$---we need to 
restrict ourselves to perfect matchings with the same number
of edges between each pair of classes. We call such perfect 
matchings \DEFN{balanced}. Observe that a complete 
$r$-partite graph $K_{r\times 2n/r}$
has a balanced perfect matching $M$ if and only if each
class is of size $(r-1)m$ for some positive integer $m$. 
In this case, $M$ has exactly $m$ edges between each 
pair of classes, and $(r-1)m=2n/r$ 
(so that $|M|=n=\binom{r}{2}m$). 

We proceed to develop a generalization of 
Theorem~\ref{3-part-thm} to the graphs 
$K_{r\times 2n/r}=K_{r\times (r-1)m}$ with constant 
$r\geq 3$ and restricting our counts to balanced 
perfect matchings. As the details are cumbersome---and don't offer further illumination---we sketch a proof
mirroring that of Theorem~\ref{3-part-thm}.
We use $\bpm(G)$ to denote the number of balanced perfect 
matchings in a graph $G$. 

There are $\binom{r}{2}$ pairs of classes in 
$K_{r \times (r-1)m}$, so following the argument that 
establishes (\ref{PM-in-tripartite}) gives the total number 
of balanced perfect matchings as
 \begin{equation}\label{PM-in-r-partite}
\bpm(K_{r \times (r-1)m}) = \left(\frac{((r-1)m)!}{(m!)^{r-1}}\right)^r (m!)^{\binom{r}{2}}.
 \end{equation}

Now suppose that $M$ is a fixed a balanced perfect matching. 
As in the proof of Theorem~\ref{3-part-thm}, we use inclusion-exclusion,
here to determine $\bpm(K_{r\times (r-1)m}-M)$.
To this end, we count the number of balanced perfect 
matchings that include edges $L\subseteq M$ as in (\ref{run-on-equation}). 
In the resulting sum (below), the indices $x_{ij}$ play 
the role of $i,j,k$ in the earlier proof; there, $(i,j,k)$ 
would be $(x_{12},x_{23},x_{13})$. If we sum over all
settings of these $\binom{r}{2}$ indices taking values 
$x_{ij}\in\{0,1,\dots,m\}$, the analogue expression of  
(\ref{run-on-equation}), which here counts 
$\bpm(K_{r \times (r-1)m} - M)$, is
\begin{equation}\label{general-run-on}
	\sum_{{\mathbf{X}_m}} (-1)^{\sum_{i < j} x_{ij}} \prod_{i < j} \binom{m}{x_{ij}}\prod_{i=1}^r \frac{((r-1)m - \sum_{j \neq i} x_{ij})!}{\prod_{j \neq i} (m-x_{ij})!}  \prod_{i < j} (m-x_{ij})! .
\end{equation}
Dividing (\ref{general-run-on}) by (\ref{PM-in-r-partite}), simplifying factors, and then applying Corollary~\ref{analysis-lemma-2} gives the following analogue of the final calculation in the proof of Theorem~\ref{3-part-thm}: 

\begin{align*}
\lim_{m \rightarrow \infty}  \frac{\bpm(K_{r \times (r-1)m}-M)}{\bpm(K_{r \times (r-1)m})} & =
\lim_{m \rightarrow \infty}\sum_{\mathbf{X}_m } (-1)^{\sum_{i < j} x_{ij}} \frac{1}{\prod_{i < j} (x_{ij})!} \left(\frac{1}{k-1}\right)^{2\sum_{j <i} x_{ij}} \\
&=
\lim_{m \rightarrow \infty} \prod_{1 \leq i < j \leq r} \sum_{0 \leq x_{ij} \leq m} \frac{(-1)^{x_{ij}}}{(x_{ij})!}  \left(\frac{1}{(r-1)^2}\right)^{x_{ij}} \\
& = \prod_{1 \leq i < j \leq r} e^{-1/(r-1)^2} \\
&= e^{-r/(2r-2)}.
\end{align*}

\bigskip

This establishes another special case of Conjecture~\ref{kayll-conj}.

\begin{theorem}
\label{k-part-thm} 
If $r\geq 3$ is fixed, and $M$ is a perfect matching in
$K_{r\times (r-1)m}=K_{r \times 2n/r}$, then \\

\hfill\(
\displaystyle{
\lim_{m \rightarrow \infty}  \frac{\bpm(K_{r \times (r-1)m}-M)}{\bpm(K_{r \times (r-1)m})} = {e^{-r/(2r-2)}}.}
\)\hfill$\Box$
\end{theorem}

\section{Cases of linear {\boldmath $r(n)$}} 
\label{linear-classes-section}

The preceding section progresses toward 
Conjecture~\ref{kayll-conj} when the number $r$ of partition classes of $K_{r\times 2n/r}$ is constant relative to 
the total number of vertices ($2n$). It's also natural 
to consider the other extreme, namely when the number of 
vertices per class is a constant, say, $c\geq 1$. In this 
case, counting vertices gives $2n=cr$, or $r=2n/c$. 
Conversely, consider any case of $r=r(n)$ being 
a linear function of $n$, say $r=Cn$. Since $r$ must divide 
$2n$, we have $C\leq 2$, and since the graphs 
$K_{r\times 2n/r}$ by definition have equal class 
sizes, we see that the number of vertices per class
is $2n/Cn$, a positive constant.
(Note, e.g., that if $C$ happens to be $1/17$, then we're
considering $34$ vertices per class, and $n\to\infty$ through 
multiples of $17$.)
In this section, we prove Conjecture~\ref{kayll-conj}
for these cases of linear $r(n)$.

We shall invoke the following identity, which follows from
results both of Godsil~\cite{Godsil81} and of
Zaslavsky~\cite{Zaslavsky81}. The statement mentions the 
\DEFN{complement $\overline{G}$} of a graph $G$, which shares $G$'s
vertex set and has the complementary edge set 
$E(K_{2n})-E(G)$; it also involves a counting function
$\mu_k(\cdot)$, which gives the number of
\DEFN{$k$-matchings} in its argument, i.e., matchings comprised of 
exactly $k$ edges. For the sake of completeness, 
we include a short proof.

\begin{lemma}[\cite{Godsil81,Zaslavsky81}]
\label{godsil-lemma}
	The number of perfect matchings in a graph $G$ on $2n$ vertices is
	\begin{equation}
	\label{godsil-zas-id}    
	\pma(G)=\sum_{k=0}^n (-1)^k \mu_k(\overline{G}) (2n-2k-1)!!.
	\end{equation}
	
\end{lemma}

\begin{proof}
   With the union of the edge sets of $G$ and $\overline{G}$ being $E(K_{2n})$, it's useful to consider the set $\MM$  of  all
   perfect matchings of  $K_{2n}$.  Using inclusion-exclusion (Theorem~\ref{PIE-stmt}), we count those $M\in\MM$
   that contain no edge of $\overline{G}$.
   So for $e\in E(\overline{G})$, let $A_e$ denote the subset of 
   $\MM$ each member of   which contains $e$. Then  PIE gives
   \[
     \pma(G)=\left|\left(\bigcup_{e\in E(\overline{G})}A_e \right)^C\right| =
\sum_{S\subseteq E(\overline{G})}(-1)^{|S|}\left|\bigcap_{e\in S} A_e\right|.
   \]
   If an index $S$ of this sum contains $k$ edges, then it contributes a nonzero term exactly when $S$ is a $k$-matching in $\overline{G}$, and there are $\mu_k(\overline{G})$  of these. Since each such $S$ spans $2k$ vertices, the remaining $2n-2k$ vertices induce a complete subgraph $K_{2n-2k}$ of
   the parent $K_{2n}$. Thus, there are $(2n-2k-1)!!$ 
   ways to complete a perfect matching starting from such an $S$,
   and this is the value of the corresponding (unsigned) term.
   Therefore, if we sum instead over the possible sizes
   $k$ of $S$, we obtain the identity (\ref{godsil-zas-id}).
\end{proof}

The expression in (\ref{kinder}) establishes the limiting proportion of perfect matchings in $K_{2n}$ with a perfect matching removed to the total number in $K_{2n}$. Our next result
generalizes this by accounting for the removal of the 
edges of a $d$-regular graph from $K_{2n}$ (a perfect matching
being the $d=1$ case). The perhaps unwieldy expression for the
vertex degrees arises because it arranges for the complementary
(`removed') graphs to be $d$-regular. We view this result
mainly as a lemma supporting Corollary~\ref{conj-case-linear-r}
below, but it may also be of independent interest.

\begin{theorem}\label{main-top-thm}
If $d$ is a fixed nonnegative integer and
$(G_{2n})$ is a sequence of $(2n-d-1)$-regular $2n$-vertex graphs,
then
\begin{equation}
\label{regular-graph-limit}
\lim_{n \rightarrow \infty} \frac{\pma(G_{2n})}{\pma(K_{2n})} = e^{-d/2}.
\end{equation}
\end{theorem}

\begin{proof}
For convenience, put $G=G_{2n}$; as noted, the complement 
$\overline{G}$ is $d$-regular. Applying Lemma~\ref{godsil-lemma} to $G$ and using $\pma(K_{2n}) = (2n-1)!!$ shows that the left fraction in
(\ref{regular-graph-limit}) is
\[
 \frac{\pma(G_{2n})}{\pma(K_{2n})} = \sum_{k=0}^n\frac{ (-1)^k \mu_k(\overline{G}) (2n-2k-1)!!}{(2n-1)!!},
\]
whose limiting value yields to
Lemma~\ref{analysis-lemma-1}. Toward that goal, we proceed to establish the lemma's three hypotheses.

The graph $\overline{G}$ contains $dn$ edges; thus,
the number $\mu_k(\overline{G})$ of its  $k$-matchings
satisfies the following (crude) upper bound:
\begin{equation}\label{mu-upper3}
\mu_k(\overline{G}) \leq \binom{dn}{k} \leq	\frac{d^k}{k!} n^k.    
\end{equation}
On the other hand, we may pick the edges of a $k$-matching one-by-one. There are $dn$ choices for the first edge $e_0$; 
since each edge is incident to $2d-2$ others, there are 
$dn-2d+1>dn-2d$ choices for the second edge (which also cannot be $e_0$). Continuing 
in this way and accounting for overlaps among the ruled-out 
edges---and for the $k!$ ways of arriving at the same $k$-matching---we find that

\begin{equation}\label{mu-lower3}
\mu_k(\overline{G}) \geq \frac{1}{k!}  \prod_{i=0}^{k-1} (dn-2d i) = \frac{d^k}{k!} \prod_{i=0}^{k-1} (n-2 i) \geq \frac{d^k}{k!}\left(n-2(k-1)\right)^k.
\end{equation}

Combining (\ref{mu-upper3}) and (\ref{mu-lower3}) gives 
\[
 \mu_k(\overline{G}) \sim \frac{d^k}{k!} n^k 
 ~~~~~~\text{(as $n\to\infty$ with $k$ constant)}.
\]
Under the same asymptotic condition, we observe that
\[
  \frac{ (2n-2k-1)!!}{(2n-1)!!}  = \prod_{i=0}^{k-1}\frac{1}{(2n-1)-2i} \sim \frac{1}{(2n)^k}.
\]
As $d$ is also constant, these asymptotic relations imply that

\begin{equation}\label{tan-app-1}
\lim_{n \rightarrow \infty} \frac{ (-1)^k \mu_k(\overline{G}) (2n-2k-1)!!}{(2n-1)!!} =\frac{(-1)^k}{k!}\left(\frac{d}{2}\right)^k.
\end{equation}
A glance at (\ref{mu-upper3}) and a little algebra gives an 
estimate of the fraction inside the limit:
\begin{equation}\label{tan-app-2}
\left| \frac{ (-1)^k \mu_k(\overline{G}) (2n-2k-1)!!}{(2n-1)!!} \right| \leq \frac{d^k}{k!}.
\end{equation}
Moreover,
\begin{equation}\label{tan-app-3}
\sum_{k=0}^\infty \frac{d^k}{k!} = e^d < \infty.
\end{equation}

With the hypotheses
(\ref{tan-app-1}), (\ref{tan-app-2}), and (\ref{tan-app-3})
established, 
we're ready to invoke Lemma~\ref{analysis-lemma-1}:
\[
\lim_{n \rightarrow \infty} \frac{\pma(G_{2n})}{\pma(K_{2n})} = \lim_{n \rightarrow \infty } \sum_{k=0}^n\frac{ (-1)^k \mu_k(\overline{G}) (2n-2k-1)!!}{(2n-1)!!} = \sum_{k=0}^\infty  \frac{(-1)^k}{k!} \left(\frac{d}{2}\right)^k = {e}^{-d/2}. \qedhere
\]
\end{proof}

\bigskip

Theorem~\ref{main-top-thm} provides just the tool we need 
to resolve Conjecture~\ref{kayll-conj} when $r(n)$ is a linear 
function of $n$. For a constant integer $c\geq 1$, let us set $r=2n/c$ as at the start of this section (recalling that $c$ 
is then the number of vertices in each partition class of 
$K_{r\times 2n/r}$). This host graph is thus 
$K_{r\times c}$, which is $(2n-c)$-regular, and if we remove a perfect matching $M$, then the resulting graph is $(2n-c-1)$-regular. Therefore, Theorem~\ref{main-top-thm} gives
\begin{align*}
\lim_{\substack{n\to\infty \\ c\mid 2n}}\frac{\pma(K_{r\times c}-M)}{\pma(K_{r\times c})} 
= \lim_{\substack{n\to\infty \\ c\mid 2n}}  \frac{\pma(K_{r\times c}-M)/\pma(K_{2n})}{\pma(K_{r\times c})/\pma(K_{2n})}= \frac{e^{-c/2}}{e^{-(c-1)/2}} = {e}^{-1/2}.
\end{align*}

\noindent
(Technically,  we twice applied Theorem~\ref{main-top-thm} to 
sequences containing the graphs in  question as subsequences.)
Thus, when $r$ is linear in $n$, we have proved the 
following partial resolution of Conjecture~\ref{kayll-conj}.

\begin{cor}
\label{conj-case-linear-r}
If $c$ is a fixed positive integer, then\\

\hfill\(
\displaystyle{
\lim_{\substack{n\to\infty \\ c\mid 2n}} \frac{\pma(K_{(2n/c)\times c}-M)}{\pma(K_{(2n/c)\times c})} = e^{-1/2}.
}\)\hfill$\Box$
\end{cor}

\section{Refinement}
\label{refinement-section}

As hinted in the preceding section,
Lemma~\ref{godsil-lemma} reformulates an integral
counting identity due to Godsil~\cite{Godsil81}.
A deeper analysis of that formula by McLeod~\cite{McLeod2009}
led her to a strong generalization of a theorem of Bollob\'{a}s~\cite{Bollobas81} from the 1980s.
One of McLeod's results takes us substantially further than Corollary~\ref{conj-case-linear-r}. A version  of  
\cite[Theorem 17]{McLeod2009} suitable for our purposes is


\begin{theorem}[\cite{McLeod2009}]
\label{mcleod-thm17}
For any given $\delta>0$, if $2 \leq d \leq O(n^{1-\delta})$
and $G$ is a $(2n -d-1)$-regular $2n$-vertex graph, then

\hfill\(\displaystyle{
\pma(G) = \frac{(2n)!}{2^{n}n!} \left( 1- \frac{d}{2n}\right)^n \cdot e^{o(1)}.
}\)\hfill$\Box$
\end{theorem}


\noindent
Notice the refinement of Theorem~\ref{main-top-thm}: the
leading fraction here is $(2n-1)!!=\pma(K_{2n})$ while
the binomial power is asymptotic to $e^{-d/2}$.

Theorem~\ref{mcleod-thm17} shows that for integers $c=c(n)$
in the range $2\leq c \leq O(n^{1-\delta})$, we have

\begin{align*}
\frac{\pma(K_{(2n/c) \times c}-M)}{\pma(K_{(2n/c)\times c})}
 =   \frac{\left( 1-c/2n\right)^n}{\left( 1-(c-1)/2n\right)^n } \cdot e^{o(1)}
 = \left(1 - \frac{1}{2n-c+1} \right) ^n \cdot e^{o(1)}.
\end{align*}
And since
\[
\left(1 - \frac{1}{2n-c+1}\right)^n\,\sim~\exp\left({-\frac{n}{2n-c+1}}\right) \,\to~ e^{-1/2}  ~~~~~~\text{(as $n\to\infty$)},
\]
our partial resolution of Conjecture~\ref{kayll-conj} is that much
closer:



\begin{cor}
\label{conj-case-frac-power-r}
If $c=c(n)$ is a positive integer with $c \leq O(n^{1-\delta})$ for a fixed $\delta>0$, then\\

\hfill\(
\displaystyle{
\lim_{\substack{n\to\infty \\ c\mid 2n}} \frac{\pma(K_{(2n/c)\times c}-M)}{\pma(K_{(2n/c)\times c})} = e^{-1/2}.
}\)\hfill$\Box$
\end{cor}

\noindent
(Of course, the case $c=1$ is already known from
Corollary~\ref{conj-case-linear-r}.)

\section{Concluding remarks}

As mentioned in the Introduction, 
Conjecture~\ref{kayll-conj} had its genesis in a heuristic  
but not rigorous argument. As such, its appeal stemmed from  
the goal of giving the Hat-Check and Kindergartner phenomena 
a common explanation. Though compelling in itself, this 
goal leaned on limited evidence. The new results presented 
here---Theorems~\ref{3-part-thm}, \ref{k-part-thm}, and 
Corollaries~\ref{conj-case-linear-r}, 
\ref{conj-case-frac-power-r}---add  substance to 
that evidence. Beyond its $r=2$ and $r=2n$ cases (the 
standard derangement and deranged matching instances),
Conjecture~\ref{kayll-conj} is now established 
for $r=3$, for balanced cases of constant $r\geq 4$,
and for $r(n)\in\Omega(n^{\delta})$ dividing $2n$ (for any fixed 
$\delta>0$).
Although this leaves open a 
gap---$\Omega(1)<r(n)<O(n^{\delta})$---we feel some satisfaction in 
having established these bookends. 

Aside from perhaps inspiring further progress toward
Conjecture~\ref{kayll-conj}, our investigations invite some 
other avenues of inquiry. We can view our parent graphs
$K_{r\times {2n}/{r}}$ as complete graphs $K_{2n}$ with
the edge sets of $r$ vertex-disjoint copies of smaller complete subgraphs $K_{2n/r}$ having been removed. So one variation on our
theme could be to study what happens if instead of these 
$E(K_{2n/r})$-removals, we remove the edge sets of triangles
($K_3$'s) or other regular structures. 
Our Theorem~\ref{main-top-thm} could provide a tool here because 
it specifically addresses graphs $G_{2n}$ obtained by the 
removal from $K_{2n}$ of the edges of a regular subgraph.
Another variation could focus the counting not on perfect 
matchings---which are $1$-regular subgraphs---rather on 
regular subgraphs of degree exceeding $1$. We have not 
attempted to establish a result analogous to
Theorem~\ref{main-top-thm} for this situation.

However, we can share a sharper estimate than one we used 
in the proof of that theorem. The following improvement of
(\ref{mu-upper3}) could conceivably be helpful in investigations
along the lines alluded above:
\[
\mu_k(\overline{G})\leq\frac{2n(2n-2)(2n-4)\cdots(2n-2(k-1))d^k}{2^{k}k!}=\binom{n}{k}d^k.
\]
We omit the proof.

Turning the focus back to perfect matching enumeration,
we should mention that the graphs $K_{r\times {2n}/{r}}$
are awfully close to the so-called 
\DEFN{Tur\'{a}n graphs}\footnote{Tur\'{a}n graphs $T_{n,r}$
feature prominently in the important subfield of extremal 
graph theory, one of the many areas championed by 
Erd\H{o}s; see, e.g., \cite{BondyMurty2008}.}
$T_{2n,r}$, which are also complete $r$-partite 
graphs---not necessarily balanced, but with the partition 
classes as close in size as possible (so pairwise all within 
a count of one). When $r$ divides $2n$---i.e., when the graphs $K_{r\times {2n}/{r}}$ are defined---the two graph classes 
coincide. But even when $r$ fails  to divide $2n$,
we believe that a variation of Conjecture~\ref{kayll-conj}
for Tur\'{a}n graphs should hold. Let it be Conjecture~2.

\vspace*{6em}

\begin{quote}
{\large Prove and conjecture!}
\hfill\textsc{Paul Erd\H{o}s} [1913--1996] {\small (attrib.)}
\end{quote}

\vfill

\pagebreak

\subsection*{Acknowledgements}
Special thanks to Mike Saks, who suggested exploring
the probabilistic heuristic, which eventually led to the formulation of Conjecture~\ref{kayll-conj}.
Emphatic thanks to Leonard Huang, who
graciously gave of his time and contributed valuable
analytic advice.
Energetic thanks to Brendan McKay, who, within half a day of this
manuscript's initial \texttt{arXiv} posting, shared the references \cite{Bollobas81},  \cite{McLeod2009}, along with suggestions how to apply them.

\bibliographystyle{abbrv}

\bibliography{JKP-biblio.bib}


\end{document}